\documentclass{tranx}
\usepackage{amssymb,latexsym, amscd}
\usepackage[all]{xy}
\usepackage{graphicx}

 \newtheorem{theorem}{Theorem}[section]
 \newtheorem{cor}[theorem]{Corollary}
 \newtheorem{lemma}[theorem]{Lemma}
 
 \theoremstyle{definition}
 \newtheorem{definition}[theorem]{Definition}
 \theoremstyle{definition}
 \newtheorem{example}[theorem]{Example}
 \theoremstyle{remark}
 
 \numberwithin{equation}{section}

\usepackage{latexsym}
\usepackage{amssymb}
\newcommand{\ben}{\begin{equation}}
\newcommand{\een}{\end{equation}}

\newcommand{\integer}{\ensuremath{{\mathbb Z}}}

\newcommand{\real}{\ensuremath{{\mathbb R}}}
\newcommand{\complex}{\ensuremath{{\mathbb C}}}

\newcommand{\rational}{\ensuremath{{\mathbb Q}}}

\newcommand{\EE}{{\mathcal E}}
\newcommand{\XX}{{\mathcal X}}

\newcommand{\KK}{{\mathcal K}}

\newcommand{\VV}{{\mathcal V}}

\newcommand{\GG}{{\mathcal G}}

\newcommand{\LL}{\mathcal{L}}
\newcommand{\MM}{\mathcal{M}}

\newcommand{\Sss}{\mathcal{S}}

\newcommand{\Loop}{\mathsf{L}}

\newcommand{\virt}{\ensuremath{{\mathrm{virt}}}}
\newcommand{\CR}{\ensuremath{{\mathrm{CR}}}}
\newcommand{\JKK}{\ensuremath{{\mathrm{JKK}}}}

\begin{document}

\title[Orbifold Cohomology and String Topology]
{Chen-Ruan Cohomology of cotangent orbifolds and Chas-Sullivan String Topology}

\author[A. Gonz\'{a}lez, E. Lupercio, C. Segovia, B. Uribe and M. A. Xicot\'{e}ncatl]{ Ana Gonz\'{a}lez, Ernesto Lupercio, Carlos Segovia, Bernardo Uribe and Miguel A. Xicot\'{e}ncatl}
\thanks{The first three and the last authors were partially supported by  Conacyt-M\'exico. The fourth author was partially supported by the ``Fondo de apoyo a investigadores j\'ovenes" from Universidad de los Andes.}

\address{Centro de Matem\'{a}tica, Universidad de la Rep\'{u}blica,
     Igu\'{a} 4225 esq. Mataojo,
    11400 Montevideo, URUGUAY}
     \address{Departamento de Matem\'{a}ticas, CINVESTAV,
     Apartado Postal 14-740
     07000 M\'{e}xico, D.F. M\'{E}XICO}
\address{Departamento de Matem\'{a}ticas, CINVESTAV,
     Apartado Postal 14-740
     07000 M\'{e}xico, D.F. M\'{E}XICO}

\address{Departamento de Matem\'{a}ticas, Universidad de los Andes,
Carrera 1 N. 18A - 10, Bogot\'a, COLOMBIA}
\address{Departamento de Matem\'{a}ticas, CINVESTAV,
     Apartado Postal 14-740
     07000 M\'{e}xico, D.F. M\'{E}XICO}

\email{ana@cmat.edu.uy, lupercio@math.cinvestav.mx, \newline
csegovia@math.cinvestav.mx, buribe@uniandes.edu.co,
xico@math.cinvestav.mx \\ }

\begin{abstract} In this paper we prove that for an almost complex orbifold, its virtual orbifold cohomology \cite{LUXSymmetric} is isomorphic as algebras to the Chen-Ruan orbifold cohomology of its cotangent orbifold.
\end{abstract}

\maketitle
\section{Introduction}

In their seminal paper Chas and Sullivan \cite{ChasSullivan} constructed a remarkable product of degree $-d$ on the homology of the free loop space $\LL M$ of an oriented manifold $M$ ,
$$ \bullet_{\mathrm{CS}} : H_*(\LL M) \otimes H_*(\LL M) \longrightarrow H_{*-d}(\LL M).$$
This product is defined using only the smooth structure of $M$ and nothing else (cf. \cite{CohenJones}).

 Viterbo \cite{Viterbo}, Salamon-Weber \cite{SW03, Weber} and Abbondandolo-Schwarz \cite{AS06} have constructed isomorphisms between a particular flavor of the Floer homology of the cotangent bundle $T^*M$ and the ordinary homology of the free loop space $$HF_*(T^*M) \simeq H_*(\LL M).$$

 Abbondandolo and Schwarz have proved that the pair of pants
product in Floer cohomology of the cotangent  corresponds to a product in the homology of the loop space, defined via
Morse theory, which Antonio Ramirez and Ralph Cohen \cite{CohenVoronov} proved is the Chas-Sullivan product. One of the main conjectures in the field  states that the symplectic field theory on the left-hand side corresponds to the string topology on the right-hand side. Here we should also mention that for a wide class of manifolds it has been shown that Floer cohomology is isomorphic to Quantum cohomology \cite{Piunikhin}.

Quite independently the study of orbifolds was revitalized by the introduction of the Chen-Ruan orbifold cohomology of a symplectic orbifold \cite{ChenRuan}.

In \cite{LupercioUribeLoopGroupoid} the second and fourth authors have
constructed a functor
$$\Loop \colon \mathbf{Orbifolds} \longrightarrow  \
S^1\mbox{-}\mathbf{Orbifolds}, \ \ \ \ \ \  \GG \mapsto \Loop \GG,$$ from
orbifolds to infinite dimensional orbifolds with actions of $S^1$.
This functor when restricted to smooth manifolds becomes the
ordinary free loop space functor $M\mapsto \LL M$, where the $S^1$
action is given by rotating the loops. But more interestingly, they have proved that the
action of $S^1$ on $\Loop \GG$ has as a fixed suborbifold $\Lambda(\GG)$
which is known as the inertia orbifold of $\GG$ (cf. \cite{LupercioUribeLoopGroupoid, deFxLuNevUri, ghost}).

 The second, the fourth and the fifth authors have proved that the \emph{homology} of the loop orbifold of a global orbifold (of the form $[Y/G]$) has the structure of a BV-algebra \cite{LUX}. The first author \cite{AnaTesis} has generalized this theorem to a general orbifold ($C^\infty$-Deligne-Mumford stack). In \cite{LUXSymmetric}, the second, the fourth and the fifth authors introduced a new version of orbifold cohomology which was called \emph{virtual} orbifold cohomology $H^*_\virt(\Lambda\GG)$. This is to be interpreted as string topology on the  ghost loop space \cite{ghost}. In \cite{LUXSymmetric} it was proved that this cohomology is a graded associative ring.  In \cite{AnaTesis} the first author has proved that it is actually a graded Frobenius algebra. In this paper we prove this result again as a corollary of our main theorem (Theorem \ref{Main}).

The main result of \cite{LUXSymmetric} stated that for an interesting family of orbifolds the Frobenius algebra $H^*_\virt(\Lambda\GG)$ is a subalgebra of the BV-algebra $H_*(\Loop \GG)$. In \cite{AnaTesis} the first author has proved that in general there is only an algebra homomorphism $$H^*_\virt(\Lambda\GG) \longrightarrow H_*(\Loop \GG)$$ with non-trivial kernel.

Chen and Ruan in their study of orbifold quantum cohomology realized that when they let the quantum deformation variable vanish they obtained a Frobenius algebra $H^*_{\mathrm{CR}}(\XX)$ for every almost symplectic orbifold $\XX$. In this paper we consider the case $\XX=T^*\GG$. The main result of this paper is:

\begin{theorem} \label{Main} For an almost complex orbifold $\GG$ we have $$H^*_{\mathrm{CR}}(T^*\GG) \cong H^*_\virt(\Lambda\GG)$$
\end{theorem}

In the case when $\GG=[Y^n/S_n]$ is the symmetric product obtained by letting the symmetric group $S_n$ act on the $n$-fold product $Y^n$ we prove that:

\begin{theorem} There is an embedding of algebras $$ H^*_\CR(T^*[Y^n/S_n]) \hookrightarrow H_*(\Loop[Y^n/S_n]).$$ \end{theorem}
The paper concludes with a conjecture regarding the Floer homology of certain Hilbert schemes.

We would like to thank conversations with A. Adem, R. Cohen, N. Ganter, Y. Ruan and C. Teleman regarding this work. We would also like to thank the MSRI, UNAM, Daniel Pineda,  Hugo Rossi and Jos\'{e} Seade for their hospitality and the terrific workshops at Berkeley, Morelia and Cuernavaca.

\section{Virtual orbifold cohomology} \label{section2}

Let $S$ be a complex manifold and let $S_1$ and $S_2$ be closed
submanifolds that intersect {\it cleanly}; that is, $U:= S_1 \cap
S_2$ is a submanifold of $S$ and at each point $x$ of $U$ the
tangent space of $U$ is the intersection of the tangent spaces of
$S_1$ and $S_2$. Let $E(S,S_1,S_2)$ be the {\it excess} bundle of
the intersection, i.e., the vector bundle over $U$ which is the
quotient of the tangent bundle of $S$ by the sum of the tangent
bundles of $S_1$ and $S_2$ restricted to $U$. Thus
$E(S,S_1,S_2)=0$ if and only if $S_1$ and $S_2$ intersect
transversally. In the Grothendieck group of vector bundles over
$U$ the excess bundle becomes
$$E(S,S_1,S_2) = T_S|_U + T_U - T_{S_1}|_U - T_{S_2}|_U.$$

Denote by $e(S,S_1,S_2)$ the Euler class of $E(S,S_1,S_2)$ and
by
\begin{eqnarray} \label{diagram excess}
\xymatrix{U \ar[r]^{i_1} \ar[rd]^h \ar[d]_{i_2} & S_1 \ar[d]^{j_1}\\
S_2 \ar[r]_{j_2} & S }
\end{eqnarray}
the relevant inclusion maps. Then for any cohomology class $\alpha
\in H^*(S_1)$ the following {\it excess intersection formula}
\cite[Prop. 3.3]{Quillen} holds in the cohomology ring of $S_2$:
\begin{eqnarray} \label{excess int formula}
j_2^* j_{1*} \alpha = i_{2*}\left(e(S,S_1,S_2)
i_1^*(\alpha)\right).\end{eqnarray}

 Consider the orbifold $[Y/G]$ where $Y$ is an almost complex manifold and $G$
  acts preserving the almost complex structure.
  Define the groups $$H^*(Y,G) := \bigoplus_{g \in G} H^*(Y^g) \times \{g\}$$
where $Y^g$ is the fixed point set of the element $g$. The group
$G$ acts in the natural way.
Denote by $Y^{g,h} = Y^g \cap Y^h$ and suppose that for every $g, h \in G$ we have cohomology classes $v(g,h) \in H^*(Y^{g,h})$, which are $G$-equivariant in the sense that $w^* v(k^{-1}gk,k^{-1}hk) = v(g,h)$ where $w : Y^{k^{-1}gk,
k^{-1}hk} \to Y^{g,h}$ takes $x$ to $w(x):=xk$. Define the map
\begin{eqnarray*}
\times : H^*(Y^g) \times H^*(Y^h) & \to & H^*(Y^{gh}) \\
(\alpha, \beta) & \mapsto & i_* \left( \alpha|_{Y^{g.h}} \cdot \beta|_{Y^{g,h}} \cdot v(g,h) \right)
\end{eqnarray*}
where $i : Y^{g,h} \to Y^{gh}$ is the natural inclusion.

Let us define now a degree shift $\sigma$ on $H^*(Y,G)$. We will declare that the degree of a class $\alpha_g \in H^*(Y^g) \subset H^*(Y,G)[\sigma]$ is $$  i + \sigma_g$$ where $$\sigma_g :=2( \dim_\complex Y - \dim_\complex Y^g),$$ and $i$ is the ordinary degree of $\alpha_g$. In this paper all dimensions and codimensions are complex. Virtual orbifold cohomology was introduced in \cite{LUXSymmetric}. There it was shown that:

\begin{theorem}
For the cohomology classes $v(g,h)=e(Y,Y^g,Y^h)$ the map
$\times$ defines an associative graded product on $H_\virt^*(Y,G):=H^*(Y,G)[\sigma]$.
\end{theorem}

\begin{definition}\label{virtual-product}
In the case when $v(g,h)=e(Y,Y^g,Y^h)$, we will call the product
$\times$ in $H^*(Y,G)$  the {\it virtual intersection product} and we will write $H^*_\virt(X,G):= (H^*(Y,G)[\sigma],\times)$.
Given that $H^*(Y,G ; \real)^G \cong H^*(I[Y/G];\real)$, the product $\times $
induces a ring structure on the orbifold cohomology of $[Y/G]$. We
will call this ring the virtual intersection ring of a global orbifold and we will denote it by $H_\virt^*(\Lambda[X/G])$.
\end{definition}

The definition of the virtual ring generalizes to a non-global orbifold. To do this we use the language of groupoids, and follow the notation of Adem-Ruan-Zhang \cite{AdemRuanZhang}. The Lemma 7.2 of \cite{AdemRuanZhang} is the generalization of the clean intersection formula of Quillen to the category of orbifolds. In the notation of \cite{AdemRuanZhang} we must replace $Y^g$ and $Y^h$ by two copies of $\Lambda \GG$, and $Y^{g,h}$ by a copy of $\GG^2$. We define  in general the virtual obstruction orbibundle $\VV\to\GG^2$ as the excess bundle of the diagram of embeddings:
\begin{eqnarray} \label{virtual diagram}
\xymatrix{\GG^2 \ar[r]^{e_1} \ar[rd]^h \ar[d]_{e_2} & \Lambda\GG \ar[d]^{j_1}\\
\Lambda\GG \ar[r]_{j_2} & \GG }
\end{eqnarray}
The definition of the degree shifting is local so we can use the same definition. We set $$H^*_\virt(\Lambda\GG):=H^*(\Lambda\GG)[\sigma].$$ The formula for the product in general becomes $$H^*_\virt(\Lambda\GG) \otimes H^*_\virt(\Lambda\GG) \longrightarrow H^*_\virt(\Lambda\GG)$$ given by $$ \alpha \times \beta := (e_{12})_*(e_1^*\alpha \cdot e_2^*\beta\cdot e(\VV)),$$
where $e_{12} : \GG^2 \to \Lambda \GG$ is the natural map that locally can be seen as the map $Y^{g,h} \to Y^{gh}$.

Ana Gonz\'{a}lez has proved in her PhD thesis the following theorem \cite{AnaTesis}.

\begin{theorem} The virtual intersection product of a general orbifold ($C^\infty$-Deligne-Mumford stack) induces the structure of a graded commutative algebra over $\rational$ on the rational cohomology $H^*_\virt(\Lambda(\GG))$, and moreover, there is a natural Frobenius algebra structure compatible with this product \end{theorem}

We will prove this theorem by different methods in this paper.

\section{Chen-Ruan Cohomology}

We will give now the definition of the Chen-Ruan cohomology following  \cite{ChenRuan}. First we need to define the degree shifting and the obstruction bundle for the Chen-Ruan theory.

Again the definition of the degree shifting is local so it is enough to define it in the case of a global quotient (cf. \cite{FantechiGotsche}).

Consider $Y$ an almost complex $G$-manifold with $G$ a finite group. Given $g\in G$ and $y \in Y^g$ we define $a(g,y)$ the \emph{age} of $g$ at $y$ as follows. Diagonalize the action of $g$ in $T_y Y$ to obtain $$g = {\mathrm{diag}}(\exp(2\pi i r_1),\ldots,\exp(2\pi i r_n)),$$ with $0\leq r_i < 1$ and set $$a(g,y):=\sum_i r_i.$$

The age $a(g,y)$ only depends on the connected component $Y^g_o$ of $Y^g$ in which $y$ lies. For this reason we can simply write $a(g,Y^g_o)$ or even $a(g)$ when there is no confusion.

Note that the age has the following interesting property $$a(g,Y^g_o)+a(g^{-1},Y^g_o) = \mathrm{codim}(Y^g_o,Y).$$

The \emph{Chen-Ruan degree shifting number} is defined then as $$s_g:=2a(g).$$ As a rational vector space the Chen-Ruan orbifold cohomology is $$H_\CR^*(Y,G):=H^*(Y,G)[s]$$ or more generally $$H_\CR^*(\GG):=H^*(\Lambda\GG)[s].$$

The definition of the obstruction bundle is modeled on the definition of the virtual fundamental class on the moduli of curves for quantum cohomology.

Let $\bar{\MM}_3(\GG)$ be the moduli space of ghost representable orbifold morphisms $f_y$ from $\mathbb{P}_3^1$ to $\GG$, where $\mathrm{im}(f)=y\in\GG_0$ and the marked orbifold Riemann surface $\mathbb{P}_3^1$ has three marked points, $z_1$, $z_2$, and $z_3$, with multiplicities $m_1$, $m_2$, and $m_3$, respectively. In \cite{AdemLeidaRuan} it is proved that $$ \bar{\MM}_3(\GG) = \GG^2.$$ Let us fix a connected component $\GG^2_o$ of $\GG^2$.

To define the Chen-Ruan obstruction bundle $\EE_o \to \GG^2_o$ we consider the elliptic complex $$\bar{\partial}_y : \Omega^0(f_y^*T\GG) \longrightarrow \Omega^{0,1}(f^*_y T\GG).$$ Chen and Ruan proved that $\mathrm{coker}(\bar{\partial}_y)$ has constant dimension along components and forms an orbivector bundle $\EE_o \to \GG^2_o$.

The formula for the Chen-Ruan product is then
$$H^*_\CR(\GG) \otimes H^*_\CR(\GG) \longrightarrow H^*_\CR(\GG)$$ given by $$ \alpha \star \beta := (e_{12})_*(e_1^*\alpha \cdot e_2^*\beta\cdot e(\EE)).$$

The following is a theorem of Chen and Ruan \cite{ChenRuan} (cf. \cite{Kaufmann}.)

\begin{theorem} $(H_\CR^*(\GG),\star)$ is a graded associative algebra, moreover it has a natural Frobenius algebra structure compatible with this product. \end{theorem}

\section{Stringy $K$-theory}

Here we should mention that both the Chen-Ruan and the virtual orbifold theories can be written in $K$-theory without much modification in the formul\ae\ \cite{JKK}. One just needs to change the Euler classes $e(\VV)$ and $e(\EE)$ for the corresponding Euler classes in $K$-theory $\lambda_{-1}(\VV)$ and $\lambda_{-1}(\EE)$ respectively. As $\integer$-modules we have $K^*_\virt(\Lambda\GG):=K^*(\Lambda\GG)$ and $K^*_\JKK(\GG):=K^*(\Lambda\GG)$.  The corresponding expressions for the products in $K$-theory are:

$$ V \times W := (e_{12})_*(e_1^*V \otimes e_2^* W \otimes \lambda_{-1}(\VV^*)),$$
and
$$ V \star W := (e_{12})_*(e_1^*V \otimes e_2^* W \otimes \lambda_{-1}(\EE^*)),$$
respectively.

\begin{theorem}[Jarvis-Kaufmann-Kimura \cite{JKK}] There exists a stringy Chern character $$\mathbf{Ch}_\JKK : K^*_\JKK(\GG) \otimes \complex \longrightarrow H^*_\CR(\GG,\complex)$$ that is a Frobenius algebra isomorphism
\end{theorem}

A similar argument shows:

\begin{theorem}[Gonz\'{a}lez \cite{AnaTesis}] There exists a virtual Chern character $$\mathbf{Ch}_\virt : K^*_\virt(\GG) \otimes \complex \longrightarrow H^*_\virt(\GG,\complex)$$ that is a Frobenius algebra isomorphism.
\end{theorem}

Note that when one has an orbifold $[Y/G]$ one can define a
$G$-Frobenius algebra $K^*_\JKK(Y,G)$ in the same way that was
defined for cohomology in section \ref{section2}  \cite{JKK}.

Define the groups $$K_\JKK^*(Y,G) := \bigoplus_{g \in G} K^*(Y^g)
\times \{g\}$$ with the natural $G$ action given by pull-back in
the first coordinate and conjugation in the second. Define the map
\begin{eqnarray*}
\star : K^*(Y^g) \times K^*(Y^h) & \to & K^*(Y^{gh}) \\
(\alpha, \beta) & \mapsto & i_* \left( \alpha|_{Y^{g.h}} \otimes
\beta|_{Y^{g,h}} \otimes \lambda_{-1}(\EE^*) \right)
\end{eqnarray*}
where $i : Y^{g,h} \to Y^{gh}$ is the natural inclusion and $\EE$
is the obstruction bundle of Chen-Ruan. Then $K^*_\JKK(Y,G)$
becomes a $G$-Frobenius algebra with the $\star$ product. Similarly we define the corresponding $K^*_\virt(Y,G)$.

\section{The cotangent orbifold}

Given an orbifold represented by a groupoid $\GG$, an orbibundle over $\GG$ is a pair $(E,\theta)$ consisting of an ordinary vector bundle over the manifold of objects $\GG_0$ and an isomorphism of the pull-backs under the source and target maps of the groupoid $\theta\colon s^* E \simeq t^* E$. To define the cotangent bundle of a real $C^\infty$-orbifold we set $E:=T^* \GG_0$, and for every arrow $g \in \GG_1$ we induce an isomorphism $$(dg^{-1})^\mathrm{T}\colon T^*_x \GG_0 \longrightarrow T^*_{gx} \GG_0.$$ All these fiber isomorphisms
assemble $\theta$.

We can always choose a groupoid $\GG$ that is both proper and \'{e}tale \cite{Moerdijk2002}. For this reason we can suppose that that action of an arrow $g\colon x \to gx$ extends to a neighborhood of $x$ as a smooth map. Moreover if $gx =x$ we can further suppose that the map $g$ acts linearly on a neighborhood of $x$ and therefore identify $d_x g$ with $g$.

In what follows we will use the notation of \cite{AdemLeidaRuan} Chapter 4 (cf. \cite{AdemRuanZhang}).

\begin{lemma} Given a real $C^\infty$-orbifold represented by a groupoid  $\GG$  then we have that $H^*(\Lambda(T^*\GG))\cong H^*(\Lambda(\GG))$ as vector spaces.
\end{lemma}
\begin{proof}
This is true because along a fixed component of a twisted sector, the inertia orbifold of the cotangent bundle is a vector bundle over the corresponding component of the twisted sector in the original orbifold. In particular, the components of $\Lambda(T^*\GG)$ are in one-to-one correspondence with those of $\Lambda(\GG)$, and are homotopy equivalent. We also note that the dimension of a component $\Lambda(T^* \GG)_o$ has exactly twice the dimension of $\Lambda(\GG)_o$. All of the above is a consequence of the orbifold isomorphism $$T^*\Lambda\GG \cong \Lambda T^* \GG.$$
\end{proof}

From now on we will assume that $\GG$ is an almost complex orbifold and we will pick a fixed (invariant) Hermitian metric on $\GG$. Such a metric induces a canonical identification $$T^*\GG \cong T\GG.$$

\begin{lemma}\label{grading} Let $s_g$ be the Chen-Ruan degree shifting number for a component of $\Lambda(T^*\GG)$ and $\sigma_g$ the virtual degree shifting number for $\Lambda(\GG)$. Then $$s_g = \sigma_g.$$
\end{lemma}
\begin{proof} This is a local statement so it is enough to show it for a global quotient $[Y/G]$. The first thing to notice is that at the zero section of the cotangent bundle we have $$T T^* Y |_Y \cong T Y \oplus T^* Y \cong TY \oplus TY.$$ Therefore the action of a group element in $G$ becomes a matrix of the form
$$\left(%
\begin{array}{cc}
  g & \  \\
  \  & g^{-1} \\
 \end{array}%
\right)
$$

From this we have
$$s_g = 2(a(g) + a(g^{-1})) = 2(\dim_\complex Y - \dim_\complex Y^g)=\sigma_g$$
as desired.
\end{proof}

\begin{theorem} Let $\KK \to \GG^2$ be the orbibundle defined by $$\KK:= T\GG^2 - e_{12}^* T\Lambda\GG$$
then we have that $$\EE = \VV + \KK$$ in $K(\GG^2)$.
\end{theorem}
\begin{proof}
Let us write the proof for  the global quotient $[Y/G]$. Let
$Z=TY\cong T^*Y$ (recall that we have picked a hermitian metric).

We will use a remarkable formula for $\EE$ obtained in \cite{JKK}
for the component $\EE^{g,h} \to Z^{g,h}$ (cf. Lemma 1.12
\cite{FantechiGotsche}).
\begin{equation}\label{jkk-formula}
\EE^{g,h} = (\Sss_g + \Sss_h - \Sss_{gh} - TZ^{gh} + TZ^{g,h})|_{Z^{g,h}}
\end{equation}
Here the bundle $\Sss_g \to Z^g$ is defined as follows. Let
$W_k^g$ be  the$k$-th eigenbundle of the action of $g$ on
$TZ|_{Z^g}$ with eigenvalue $0 \leq r_k <1$. Then we define
$\Sss_g$ as the sum over all eigenbundles: $$\Sss_g:= \bigoplus_k
r_k W_k^g.$$ Formula \ref{jkk-formula} is a non-trivial
consequence of the Eichler trace formula. Similarly we define $S_g
\to Y^g$.

We need to compute $\EE^{g,h}|_{Y^{g,h}}$. For this we notice
first that by equation (28) of \cite{JKK} we have: $$\Sss_g|_{Y^g}
= S_g + S_{g^{-1}} = N_g = (TY - TY^g)|_{Y^g}.$$ Similarly we can
write $$ TZ^{gh}|_{Y^{g,h}} = 2TY^{gh}|_{Y^{g,h}}$$ and $$
TZ^{g,h}|_{Y^{g,h}} = 2TY^{g,h}.$$ Putting all this back together
in equation                  \ref{jkk-formula} we get
$$\EE^{g,h}|_{Y^{g,h}} = (TY - TY^g + TY - TY^h  - TY + TY^{gh} -
2TY^{gh} + 2TY^{g,h})|_{Y^{g,h}} =$$ $$ = (TY -
TY^g-TY^h+TY^{g,h})|_{Y^{g,h}} + (TY^{g,h}-TY^{gh})|_{Y^{g,h}} =
$$ $$ = \VV^{g,h} + \KK^{g,h}.$$ The proof in the general case is
identical replacing $Y^g$ and $Y^h$ by $\Lambda\GG$ and $Y^{g,h}$
by $\GG^2$. We leave the details to the reader.

\end{proof}

With this result and Lemma 7.2 of \cite{AdemRuanZhang} we can prove the following theorem.

\begin{theorem} The zero-section inclusion $j:\Lambda \GG \to T^* \Lambda \GG$ induces
an isomorphism of rings $$ j^* : K^*_\JKK(T^* \Lambda \GG) \longrightarrow   K^*_\virt(\Lambda\GG)
 .$$ Similarly the zero
section $j:Y\to T^*Y$ induces a ring isomorphism $$ j^* : K^*_\JKK(T^*Y,G) \longrightarrow
K^*_\virt(Y,G)  .$$
\end{theorem}

\begin{proof}
We are using again the notations from \cite{AdemRuanZhang}.
Consider the following commutative diagram of natural inclusions:
\begin{eqnarray} \label{main-diagram}
\xymatrix{
\Lambda\GG \times \Lambda\GG \ar[d]^{j\times j} &  \GG^2 \ar[l]_{\ \ \  \ \ \ \Delta} \ar[d]^{j'} \ar[r]^{\iota} & \Lambda\GG \ar[d]^{j}\\
T^*\Lambda\GG \times T^*\Lambda\GG &  T^*\GG^2 \ar[l]_{\ \ \ \ \ \
\ \Delta^\tau} \ar[r]^{\iota^\tau} & T^*\Lambda\GG, }
\end{eqnarray}
where $\Delta^\tau:=de_1 \times de_2$, $\Delta:=e_1\times
e_2$,$\iota:=e_{12}$, and $\iota^\tau = de_{12}$.

The excess intersection bundle of the right hand side square is:
$$TT^*\Lambda\GG + T\GG^2 - T\Lambda\GG-TT^*\GG^2 =  2T\Lambda\GG + T\GG^2 - T\Lambda\GG-2T\GG^2=-\KK.$$

Then we have:
$$j^*(\alpha \star \beta) = j^*\iota^\tau_*({\Delta^\tau}^*(\alpha,\beta)\otimes \lambda_{-1}(\EE^*)) =
 \iota_* {j'}^* ({\Delta^\tau}^*(\alpha,\beta) \otimes \lambda_{-1}(\EE^*) \otimes \lambda_{-1}(-\KK^*)) =$$

$$=
\iota_* (\Delta^*(j^* \alpha,j^* \beta)\otimes
\lambda_{-1}({j'}^*\VV^*)) = j^*(\alpha) \times j^*(\beta).$$ This
proves the first part of the theorem. For the second part of the
theorem we just need to consider the following diagram:

\begin{eqnarray}
\xymatrix{
Y^g \times Y^h \ar[d]^{j\times j} &  Y^{g,h} \ar[l]_{\ \ \ \ \ \ \Delta} \ar[d]^{j'} \ar[r]^{\iota} & Y^{gh} \ar[d]^{j}\\
TY^g \times TY^h &  TY^{g,h} \ar[l]_{\ \ \ \ \ \ \ \Delta^\tau}
\ar[r]^{\iota^\tau} & TY^{gh}. }\end{eqnarray} We leave the
details to the reader.
\end{proof}

Writing this proof again in cohomology, or using Lemma
\ref{grading} and the Chern character, we obtain Theorem
\ref{Main}.

\section{Final remarks}

We have several immediate consequences of Theorem \ref{Main}. The following one is a property that
$H^*_{\virt}(\Lambda\GG)$ inherits from the Chen-Ruan theory.

\begin{cor} $H^*_{\virt}(\Lambda\GG)$ has the structure of a graded Frobenius algebra. In the case
$\GG = [Y/G]$ we have that $H^*_\virt(Y,G)$ has the structure of graded a $G$-Frobenius algebra. \end{cor}
\begin{proof} This is a consequence of the results of Kaufmann \cite{Kaufmann} and of Jarvis-Kaufmann-Kimura \cite{JKK} for Chen-Ruan theory. \end{proof}

Ana Gonz\'{a}lez has proved this theorem directly in her PhD dissertation \cite{AnaTesis}.


The following two examples follow from the corresponding statements for the Chen-Ruan theory.

\begin{example}
 When $G={1}$ and the orbifold is actually a manifold, the virtual intersection ring coincides with the usual intersection ring of a smooth manifold, which is a Frobenius algebra. Here the statement of Theorem \ref{Main} states that $$H^*(T^*M) \cong H^*(M)$$ which is true by the homotopy invariance of cohomology.
\end{example}
\begin{example}
When $Y=\{\bullet\}$ is a point, the virtual intersection ring becomes the Dijkgraaf-Witten  Frobenius algebra associated to a finite group \cite{DijkgraafWitten}. Here the cotangent orbifold equals the original orbifold and the statement of Theorem \ref{Main} is a tautology.
\end{example}

\begin{example} This is a more interesting example. Consider the global orbifold
$$\GG=[Y^n/S_n],$$ namely the $n$-th symmetric product of a complex manifold.

A theorem of Fantechi-G\"{o}ttsche \cite{FantechiGotsche} and Uribe \cite{Uribe} states that when $X$ is a complex projective surface with trivial canonical class the Chen-Ruan cohomology of the symmetric product orbifold is isomorphic to the ordinary cohomology of the Hilbert scheme (which is a resolution of singularities) $X^{[n]}$, namely $$H^*_{\CR}([X^n/S_n]) \cong H^*(X^{[n]}).$$

The main theorem of \cite{LUXSymmetric} states that there is an embedding of rings from virtual orbifold cohomology to orbifold string topology of the loop orbifold $\Loop[Y^n/S_n]$, $$H^*_\virt(\Lambda[Y^n/S_n]) \hookrightarrow H_*(\Loop[Y^n/S_n])$$
where the right hand side carries the product defined in \cite{LUX}.

Combining these results with theorem \ref{Main} we obtain:

\begin{theorem} There is an embedding of algebras $$ H^*_\CR(T^*[Y^n/S_n]) \hookrightarrow H_*(\Loop[Y^n/S_n]).$$ \end{theorem}

We end this paper with a conjecture. We believe that when $Y$ is a complex curve there is an isomorphism
$$ FH_*((T^*Y)^{[n]}) \cong H_*(\Loop[Y^n/S_n]).$$
between the Floer homology of the Hilbert scheme of $T^*Y$  and
the orbifold string topology of the loop symmetric product. For a
general complex manifold $Y$ it is reasonable to believe that
$$ {FH_*}^\CR(T^*[Y^n/S_n]) \cong H_*(\Loop[Y^n/S_n]),$$
where the left hand side is the yet to be defined Chen-Ruan
version of Floer homology.

\end{example}

\bibliographystyle{amsplain}
\bibliography{Loop_orb}

\end{document}